\documentclass[letterpaper, 10 pt, conference]{ieeeconf}
\IEEEoverridecommandlockouts                              
\usepackage{amsmath,mathrsfs,amssymb,float,color}
\usepackage{graphicx}
\usepackage{verbatim}
\usepackage{ upgreek }
\usepackage[T1]{fontenc}
\usepackage[utf8]{inputenc}
\usepackage{epsfig,subfigure,fancybox,balance}
\usepackage{diagbox}
\usepackage{youngtab}

\def\para#1{\vskip 0.125\baselineskip\noindent{\it #1}}
\def\qed{\strut\hfill $\Box$}

\newcommand{\beq}[1]{\begin{equation} \label{#1}}
\newcommand{\eeq}{\end{equation}}
\newcommand{\bea}{\bed\begin{array}{rl}}
\newcommand{\eea}{\end{array}\eed}
\newcommand{\bed}{\begin{displaymath}}
\newcommand{\eed}{\end{displaymath}}
\newcommand{\barray}{\begin{array}{ll}}
\newcommand{\earray}{\end{array}}
\newcommand{\disp}{\displaystyle}
\newcommand{\ad}{&\!\disp}
\newcommand{\aad}{&\disp}
\newcommand{\al}{\alpha}
\newcommand{\e}{\varepsilon}

\newcommand{\la}{\lambda}
\newcommand{\La}{\Lambda}
\newcommand{\sg}{\sigma}

\newcommand{\ga}{\gamma}

\newcommand{\Ga}{\Gamma}
\newcommand{\dl}{\delta}
\newcommand{\Dl}{\Delta}
\newcommand{\cd}{(\cdot)}

\def\phi{\varphi}
\def\indi{{\bf 1}}

\newcommand{\CA}{{\mathcal A}}
\newcommand{\CF}{{\mathcal F}}
\newcommand{\CM}{\mathcal{M}}

\newcommand{\CX}{\mathcal{X}}
\newcommand{\CC}{\mathcal{C}}
\newcommand{\CO}{\mathcal{O}}

\newcommand{\CL}{\mathcal{L}}
\newcommand{\CP}{\mathcal{P}}
\newcommand{\CE}{\mathcal{E}}
\newcommand{\CK}{\mathcal{K}}

\newcommand{\CS}{\mathcal{S}}
\newcommand{\CV}{\mathcal{V}}

\newcommand{\CH}{\mathcal{H}}
\newcommand{\CN}{\mathcal{N}}

\newcommand{\CG}{\mathcal{G}}
\newcommand{\EE}{{\mathbb E}}
\newcommand{\PP}{{\mathbb P}}

\newcommand{\NN}{{\mathbb N}}

\newcommand{\rr}{{\mathbb R}}

\newcommand{\wdt}{\widetilde}
\newcommand{\wdh}{\widehat}

\newtheorem{thm}{\rm{\bf{Theorem}}}
\newtheorem{lem}[thm]{\rm{\bf{Lemma}}}
\newtheorem{defn}[thm]{\rm{\bf{Definition}}}

\newtheorem{prop}[thm]{\rm{\bf{Proposition}}}
\newtheorem{rem}[thm]{Remark}

\newtheorem{const}[thm]{\rm{\bf{Construction}}}

\newcommand{\thmref}[1]{Theorem~{\rm \ref{#1}}}
\newcommand{\lemref}[1]{Lemma~{\rm \ref{#1}}}

\newcommand{\propref}[1]{Proposition~{\rm \ref{#1}}}

\newcommand{\secref}[1]{Section~{\rm \ref{#1}}}
\newcommand{\constref}[1]{Construction~{\rm \ref{#1}}}

\title{\LARGE  \bf Large Deviations Analysis For Regret Minimizing Stochastic Approximation Algorithms
}
\author{Hongjiang Qian and Vikram Krishnamurthy
\thanks{This work was supported by U. S. Army Research Office  grant W911NF-24-1-0083 and NSF grant CCF-2112457.} 
\thanks{Hongjiang Qian {\tt\small{hongjiang.qian@cornell.edu}} and Vikram Krishnamurthy {\tt\small{vikramk@cornell.edu}} are with the School of Electrical and Computer Engineering, Cornell University, Ithaca, NY 14853, USA.}
}

\begin{document}
\maketitle
\thispagestyle{empty}
\pagestyle{empty}

\begin{abstract}
Motivated by learning of correlated equilibria in non-cooperative games, we perform a large deviations analysis of a regret minimizing stochastic approximation algorithm. The regret minimization algorithm we consider comprises multiple agents that communicate over a graph to coordinate their decisions. We derive an exponential decay rate towards the algorithm's stable point using large deviations theory. Our analysis leverages the variational representation of the Laplace functionals and weak convergence methods to characterize the exponential decay rate.
\end{abstract}

\section{Introduction}
This paper is motivated by regret minimization algorithms for 
 finite repeated non-cooperative games. Such games are parameterized by 
\beq{game}
\CG=\{\CK, (\CA^k)_{k\in \CK}, \CC, (U^k)_\CK, (\boldsymbol{\sg}^k)_{k\in \CK}\},
\eeq
where $\CK=\{1,2,\dots, \bar K\}$ is the set of agents and $k \in \CK$ represents individual agent. The $\CA^k=\{1,2,\dots, A^k\}$ denotes the set of actions for each agent $k$ and $|\CA^k|=A^k$. We use $a^k\in \CA^k$ to denote the generic action of agent $k$. The $\CC=(\CV,\CE)$ is an undirected graph, where $\CV=\CK$ is the set of agents. Define $\CN_k=\{k' \in \CV \setminus \{k\}; (k,k')\in \CE\}$ and $\CN_k^c=\CN_k \cup \{k\}$, the open and closed neighborhoods for each agent $k$. Denote by $\CS_k=\CK \setminus \CN_k^c$ the set of non-neighbors of agent $k$. Each agent $k$ exchange information only with the neighbors $\CN_k$ that are determined by the connectivity graph $\CC$. We assume the  neighborhood monitoring assumptions:
$(k,k')\in \CE \Leftrightarrow k$  knows $ a_n^{k'}$  and  $k'$  knows $ a_n^k$ at the end of period  $n$.

The $U^k: \CA^\CK \to \rr$ is the bounded payoff function for agent $k$ and $\CA^{\CK}$ is the set of $\bar K$-tuples of action profiles, where a generic element of $\CA^\CK$ is denoted by $\boldsymbol{a}=(a^1,\dots, a^{\bar K})$. Following the common notation in game theory, we can rewrite $\boldsymbol{a}$ as $(a^k, \boldsymbol{a}^{-k})$, where $\boldsymbol{a}^{-k}\in \times_{k'\in \CK\setminus \{k\}}\CA^{k'}$. The $\boldsymbol{a}^{-k}$ can be further rewrite as $(\boldsymbol{a}^{\CN_k}, \boldsymbol{a}^{\CS_k})$, where $\boldsymbol{a}^{\CN_k}\in \CA^{\CN_k}=\times_{k'\in \CN^k}\CA^{k'}$ is the action profile of agent $k$'s neighbors and $\boldsymbol{a}^{\CS_k}\in \CA^{\CS_k}=\times_{k'\in \CS_k} \CA^{k'}$ is the action profile of agent $k$'s non-neighbors. The payoff function $U^k$ has two parts: $U_l^k(a^k, \boldsymbol{a}^{\CN_k})$ due to the local interaction with neighbors and $U_g^k(a^k, \boldsymbol{a}^{\CS_k})$ due to strategic global interaction with non-neighbors. That is,
\beq{U}
U^k(a^k,\boldsymbol{a}^{-k})=U_l^k(a^k,\boldsymbol{a}^{\CN_k})+U_g^k(a^k,\boldsymbol{a}^{\CS_k}),
\eeq
where $U_l^k(a^k,\boldsymbol{a}^{\CN_k})=0$ if $\CN_k=\emptyset$. Here, we suppose agent $k$  knows the functions $U_l^k\cd$ and $U_g^k\cd$.

The $\boldsymbol{\sg}^k$ is randomized strategies selected by agent $k$. It belongs to the simplex $\Dl \CA^k=\{\boldsymbol{p}^k\in \rr^{A^k}: \boldsymbol{p}^k(a)\geq 0, \sum_{a\in \CA^k} \boldsymbol{p}^k(a)=1 \}$ and is a map
$
\boldsymbol{\sg}_{n+1}^k : \bigcup_n \CH_n^k \to \Dl \CA^k,
$
where $\CH_n^k=(\times_{k'\in \CN_k^c} \CA^{k'}, U^k)^n$ denotes the space of all possible joint moves of agent $k$, his neighbors $\CN_k$, and possible payoffs to agent $k$ up to period $n$. Denote by $\CF_n^k=\sg\{h_n^k : h_n^k\in \CH_n^k\}$. We use $\al_n^k(i,j), \beta_n^k(i,j)$ to represent the weighted time-averaged local-regret and weighted time-averaged global-regret had the agent selected action $j$ every time he played action $i$ in the past, respectively. That is,
\beq{regret}\barray
\al_n^k(i,j)\ad =\e \sum_{\iota=1}^n (1-\e)^{n-\iota} \\
\aad\qquad  \times [U_l^k(j, \boldsymbol{a}_\iota^{\CN_k})-U_l^k(a_\iota^k, \boldsymbol{a}_\iota^{\CN_k})] \indi_{\{a_\iota^k=i\}},\\
\beta_n^k(i,j) \ad =\e \sum_{\iota=1}^n (1-\e)^{n-\iota} \\
\aad \qquad \times [U_g^k(j,\boldsymbol{a}_\iota^{\CS_k}) - U_g^k(a_\iota^k, \boldsymbol{a}_\iota^{\CS_k})] \indi_{\{a_\iota^k = i\}}.
\earray\eeq
Here $0<\e \ll 1$ is called adaptation rate. The transition probability of $a_n^k$ is defined by $\psi_n^k(i):=\PP(a_n^k=i|a_{n-1}^k)=(1-\kappa)\min\{|\al_n^k(a_{n-1}^k, i)+\beta_n^k(a_{n-1}^k,i)|^+/\xi^k, 1/A^k\}+\kappa/A^k$ if $i\neq a_{n-1}^k$, and  $\PP(a_n^k=i|a_{n-1}^k)=1-\sum_{j\in \CA^k\setminus \{i\}}\psi_n^k(j)$ if $i=a_{n-1}^k$. The $\kappa \in (0,1)$ is the exploration factor and $\xi^k> A^k \cdot |U^k_{\text{max}}-U^k_{\text{min}}|$ is a constant, where $U_{\text{max}}^k$ and $U_{\text{min}}^k$ denote the upper and lower bounds on the payoff function for agent $k$, respectively.
It can be seen $\al_n^k(i,j)$ and $\beta_n^k(i,j)$ are updated according to the following regret-based  stochastic approximation algorithm:
\beq{sa-game}\barray
\!\!\!\!\!\!\! \al_n^k(i,j)\ad \!\!\!\!=\al_{n-1}^k(i,j)+ \e ([U_l^k(j,\boldsymbol{a}_n^{\CN_k}) - U_l^k(a_n^k, \boldsymbol{a}_n^{\CN_k})]\\
\aad\qquad\qquad\qquad\qquad \times \indi_{\{a_n^k=i\}}-\al_{n-1}^k(i,j)), \\
\!\!\!\!\!\!\! \beta_n^k(i,j)\ad\!\!\!\! = \beta_{n-1}^k(i,j)+ \e ([U_g^k(j,\boldsymbol{a}_n^{\CS_k}) - U_g^k(a_n^k, \boldsymbol{a}_n^{\CS_k})] \\
\aad \qquad\qquad\qquad\qquad\times \indi_{\{a_n^k=i\}} - \beta_{n-1}^k(i,j)).
\earray\eeq
Define the piecewise constant interpolation $(\al^{k,\e}, \beta^{k,\e})$ by $(\al^{k,\e}(t),\beta^{k,\e}(t)):=(\al_n^k, \beta_n^k)$ if $t\in [n\e, (n+1)\e), t\in [0,T]$. In \cite{NKY13}, Namvar et al.  proved  that $(\al^{k,\e}\cd, \beta^{k,\e}\cd )$ converges weakly to $(\al^k\cd, \beta^k\cd)$ satisfying the differential inclusions:
\beq{sa-limit}\left\{\barray
\disp \frac{d \al^k}{dt} \in \CL_1^k(\al^k, \beta^k)-\al^k, \\
\disp \frac{d \beta^k}{dt} \in  \CL_2^k(\al^k, \beta^k)-\beta^k,
\earray
\right.\eeq
with suitable operators $\CL_1^k$ and $\CL_2^k$. The appearance of differential inclusion in \eqref{sa-limit} is because agents are oblivious to strategies where actions of neighbors and non-neighbors are drawn. In fact, the result in \cite{NKY13} is more general in that instead of considering a static game, they studied a dynamic game where the the game parameters evolve according to a finite-state Markov chain and the limiting system \eqref{sa-limit} becomes interconnected differential inclusions; see \cite{NKY13} and references therein for details. Once we have the convergence result, a natural question is how fast the pair $(\al^{k,\e}, \beta^{k,\e})$ converges to the stationary point $(\al^*, \beta^*)$ of \eqref{sa-limit}. Such question can be answered by the classical rate of convergence (central limit type) using the asymptotic normality of $((\al^{k,\e}, \beta^{k,\e})-(\al^*,\beta^*))/\sqrt{\e}$; see \cite{KY03} and \cite{YKI04}. This asymptotic result is useful, but does not provide enough information on the convergence rate. Once $(\al_n^\e, \beta_n^\e)$ is in a small fixed neighborhood $G$ of $(\al^*,\beta^*)$, the asymptotic normality is not as important as the probability escape from $G$. This escape probability could provide a better measurement of the stability and quality of algorithms.

In this paper, we aim to study the useful non-classic estimates of ``rate of convergence'' for $(\al_n^k, \beta_n^k)$, via the theory of large deviations.
Instead of studying the general dynamic game in \cite{NKY13}, we focus on the static game. That is, the game parameters do not evolve. Moreover, for the moment, we only consider the limiting systems \eqref{sa-limit} are differential equations. In other words, we assume the probability distribution of the actions of neighbors and non-neighbors are known to us. The analysis when \eqref{sa-limit} are differential inclusions is more involved and we leave it for future work. 

To proceed, we consider the following general algorithm:
\beq{sa}
\left\{\barray
X_{n+1}^1 \ad =X_n^1+ \e (U_1(\Phi_n, \Psi_{n}^1)-X_n^1),\\
X_{n+1}^2 \ad =X_n^2+ \e (U_2(\Phi_n,\Psi_{n}^2)-X_n^2),
\earray\right.
\eeq
where $X_n^1, X_n^2\in \rr^r, r=A^k\times A^k$. The algorithm \eqref{sa-game} can be considered as a special case of \eqref{sa}, where $X_n^1, X_n^2$ represent $\al_n^k$ and $\beta_n^k$, respectively. The $\Phi_n$ denotes agent $k$'s action $a_n^k$, $\Psi_n^1$ and $\Psi_n^2$ represent the joint action profile of his neighbors $\boldsymbol{a}_n^{\CN_k}$ and non-neighbors $\boldsymbol{a}_n^{\CS_k}$, respectively. 

Throughout the paper, we assume $\{\Phi_n\}$, $\{\Psi_n^1\}$, and $\{\Psi_n^2\}$ take values in the finite-state space $\CM=\{1,\dots, m\}$. The $\{\Psi_n^1\}$ and $\{\Psi_n^2\}$ are independent of $\{\Phi_n\}$. Let $X_n=(X_n^1, X_n^2)'$ and $\Psi_n=(\Psi_n^1,\Psi_n^2)'$, where $z'$ denotes the transpose of $z$. We define
$
U(X_n, \Phi_n, \Psi_n):=(\wdh U_1(X_n^1,\Phi_n, \Psi_n^1), \wdh U_2(X_n^2, \Phi_n, \Psi_n^2))',
$
where
\beq{hat-U}\barray
\wdh U_1(X_n^1, \Phi_n, \Psi_n^1) \ad := U_1(\Phi_n, \Psi_n^1)-X_n^1, \\
\wdh U_2(X_n^2, \Phi_n, \Psi_n^2) \ad := U_2(\Phi_n, \Psi_n^2)-X_n^2.
\earray\eeq
Rewriting \eqref{sa} further yields the following concise form:
\beq{sa-z}
X_{n+1}=X_n+ \e U(X_n, \Phi_n, \Psi_n),
\eeq
where $X_n\in \rr^d, d=2r$. Suppose $\{X_n, \Phi_n\}$ is a Markov chain with transition probability given by
\beq{rho-1}
\PP(\Phi_{n+1}=j|\Phi_n =i , X_n=x)= \rho^1_x(i,j). 
\eeq
Such $\{\Phi_n\}$ is known as state-dependent noise. Moreover, we assume $\{\Psi_n\}$ is a Markov chain with state space $\CM^2:=\CM\times \CM$ and with transition probability $\rho^2$ that is independent of $\{\Phi_n\}$. For the rest of the paper, we will work with \eqref{sa-z}.

{\bf Context}.
The study of the asymptotic behavior of stochastic approximation algorithms via large deviations can be traced back to the work of Kushner \cite{Ku84-tac}. In his study, Kushner used a decreasing step size $\e_n=(n+1)^{-\rho}$, where $0 < \rho\leq 1$, and considered the function $U(X_n,\Phi_n,\Psi_n)=\wdt U(X_n)+\Phi_n$ for some function $\wdt U$ with ${\Phi_n}$ being i.i.d. zero-mean Gaussian noise. Subsequently, Dupuis and Kushner extended these findings to scenarios involving  non-additive noise in \cite{DK85}, and to projected or constrained  stochastic approximation algorithm on a bounded convex set in \cite{DK87, Dup87}. Further investigations into the rate of convergence through large deviations for projected stochastic approximation algorithm with $U(X_n, \Phi_n, \Psi_n)=U_n(X_n,\Phi_n)$ was conducted in \cite{DK89}, where $\{U_n\}$ is a sequence of vector-valued mutually independent, but not necessary stationary, random fields parameterized by $X_n, \Phi_n$. Moreover, in \cite{Du88}, Dupuis investigated the large deviations for stochastic approximation algorithm with state-dependent noise and discontinuous forcing term $U(X_n,\Phi_n,\Psi_n)=U(X_n,\Phi_n)$. The underlying idea in these findings was the identification of a Hamiltonian $H(x,\al)$ as a limiting log-moment generation function, with the local rate function $L(x,\beta)$ defined as the convex conjugate of $H(x,\al)$. We refer to \cite{DZ09} and \cite{DE97} for comprehensive study of large deviations.

Compared with existing literature, our stochastic approximation algorithm \eqref{sa} or \eqref{sa-z} incorporates two independent  non-additive noise $\{\Phi_n\}$ and $\{\Psi_n\}$, with $\{\Phi_n\}$ being state-dependent while $\{\Psi_n\}$ is not. This setting is natural in the context of regret minimization \eqref{sa-game} in the finite repeated non-cooperative games. We mention that, in general, the state spaces of $\{\Phi_n\}$ and $\{\Psi_n\}$ are not necessarily compact (in our case the finite-state space $\CM$). The analysis for non-compact state space can also be employed in a similar way; see \cite{DE97}. Since our motivation comes from learning in non-cooperative games, we use the finite-state space $\CM$ as the state spaces of $\{\Phi_n\}$ and $\{\Psi_n\}$.

Rather than pursuing the identification of a Hamiltonian, we reformulate the large deviations problem as a stochastic control problem and employ the weak convergence method for its analysis, following the approach in \cite{DE97}. The weak convergence method surpasses the Freidlin–Wentzell approach \cite{FW12} by eliminating the need for establishing exponential probability estimates, a more complex task, and by allowing for the explicit expression of the local rate function in terms of the transition kernels $\rho_x^1$ and $\rho^2$, which is normally implicit in Freidlin–Wentzell approach. This method also has advantage over the approach in Feng and Kurtz \cite{FK06}, where it is necessary to establish the exponential tightness and the uniqueness of an infinite dimensional Hamilton-Jacobi equations. The weak convergence method used in this paper only requires the establishment of tightness.

The paper is organized as follows. In \secref{sec:main}, the definition of large deviations principle and hypothesis are stated. The main results are also presented in this section. In \secref{sec:up-bound}, the large deviations upper bound is proved via weak convergence methods. The \secref{sec:low-bound} is devoted to proving the lower bound of large deviations.

\section{Preliminaries and main results}\label{sec:main}
In this paper, $C_{x_0}([0,T];\rr^d)$ denotes the space of continuous function from $[0,T]$ to $\rr^d$ with initial value $x_0$, equipped with the uniform topology. Denote by $\|\cdot\|$ the norm in $\rr^d$. The $\CP(\CM)$ denotes the space of probability measures on $\CM$. Given a probability function $r$, we denote, abusing notation, $\sum_{\CM} f(x)r(x)$ by $\int_{\CM} f(x)r(dx)$ for a real function $f$ on $\CM$. For probability measures $r_1$ on $\CM $ and $r_2$ on $\CM^2$, the product measure on $\CM\times \CM^2$ is defined by $(r_1\otimes r_2)(dz_1,dz_2):=r_1(dz_1)\otimes r_2(dz_2)$. We omit the symbol $\otimes$ if there is no confusion. We use $A, B, C, D, E$ to denote the generic measurable set in $\CM, \CM, \CM^2, \CM^2, [0,T]$, respectively. For simplicity, we let $N:=T/\e$ be an integer and set $t_i=i \e, i=0,\dots, N$. We also use $t/\e$ to denote the integer $\lfloor t/\e \rfloor$ for notation simplicity.
The $K$ represents a generic constant that changes from place to place.

For a Polish space $\CX$, a function $I:\CX \to [0,\infty]$ is said to be a rate function on $\CX$ if for each $c<\infty$, the level set $\{x\in \CX: I(x)\leq c\}$ is a compact subset of $\CX$. Let us recall the following Laplace principle, which is equivalent to the large deviations principle; see \cite{DE97}.
\begin{defn}\label{defn:laplace}
Let $I$ be a rate function on $\CX$. A sequence $\{X^\e\}$ is said to satisfy the Laplace principle on space $\CX$ with speed $\e^{-1}$ and rate function $I$, if for all bounded continuous function $F: \CX \to \rr$, we have the Laplace principle upper bound
\bea\ad 
\limsup_{\e\to 0} \e \log\EE e^{- \e^{-1} F(X^\e)}\leq  -\inf_{x\in \CX} \big\{F(x)+I(x)\big\},
\eea
and the Laplace principle lower bound
\bea \ad
\liminf_{\e\to 0} \e \log\EE e^{- \e^{-1} F(X^\e)}\geq  -\inf_{x \in \CX} \big\{F(x)+I(x)\big\}.
\eea
\end{defn}
To proceed, we make the following assumptions:

\para{\rm{(A1)}}\ For each $z_1\in \CM, z_2\in \CM^2$, the mapping $x\mapsto U(x, z_1, z_2)$ is Lipschitz continuous, uniformly with respect to $(z_1,z_2)$ and there exists a $M>0$ such that 
\bea\ad 
\sup_{(x,z_1,z_2)\in \rr^d\times \CM \times \CM^2} |U(x, z_1, z_2)| \leq M <\infty. 
\eea

\para{\rm{(A2)}}\  For each $z_1 \in \CM$, the mapping $x \mapsto \rho_{x}^1(z_1,\cdot)$ from $\rr^{d}$ to $\CP(\CM)$ is continuous. For each $x\in \rr^d$, $\rho_x^1$ and $\rho^2$ are irreducible and aperiodic.

\begin{rem}
We mention that $\CM$ is a compact space, then $\CP(\CM)$ is compact as is $\CP(\CM^3)$. The condition (A2) implies for each $x\in \rr^d$, $\rho_x^1$ and $\rho^2$ have a unique invariant measure denoted by $\pi_x^1$ and $\pi^2$, respectively. In addition, the transitivity condition is satisfied for the chain $(\Phi_n,\Psi_n)$: for any $x\in \rr^d$, there exists positive integer $\ell_0$ and $n_0$ such that for all $(z_1,z_2)\in \CM^3$ and $(z_3,z_4)\in \CM^3$, 
\beq{tr-cond}
\sum_{i=\ell_0}^{\infty}\frac{1}{2^i} \rho_x^1(z_1,d\theta)\rho^2(z_2,d\psi) \leq\!\! \sum_{j=n_0}^\infty \frac{1}{2^j} \rho_x^1(z_3,d\theta)\rho^2(z_4,d\psi).
\eeq
\end{rem}

Define the continuous time piecewise constant interpolation $X^\e(t)$ by $X^\e(t):=X_n^\e$ if $t\in [n\e, (n+1)\e)$ with $X_0^\e=x_0$. Then under assumptions (A1) and (A2), as in \cite{NKY13}, $X^\e$ converges weakly to $X$, where $
\dot{X}(t)=\bar U(X(t))$ and $\bar U(x):=\int_{\CM^3} U(x,z_1,z_2)\pi_x^1(dz_1)\pi^2(dz_2)$ for $x\in \rr^d$.
We now introduce the following construction, which is useful to identify quantities for the representation formula in Laplace principle and the weak convergence analysis.

\begin{const}\label{const}\rm{
Suppose we are given a probability measure $\mu^\e\in \CP((\CM^3)^N)$  and decompose it in terms of conditional distribution $[ \mu^\e]_{i|0,1,\dots,i-1}$ on the $i$-th variables given variables $0$ to $i-1$ as following:
\bea\ad\!\!\!\!\!
\mu^\e(dq_0, dq_1, \cdots, dq_{N-1}) \!=[\mu^\e]_0(dq_0)[\mu^\e]_{1|0}(dq_1|q_0) \otimes \cdots \\
\aad\qquad\qquad\quad \otimes [\mu^\e]_{N-1|1,\dots, N-2}(dq_{N-1}|q_0,\dots, q_{N-2}).
\eea
Let $(\bar q_0^\e,\bar q_1^\e,\dots, \bar q_{N-1}^\e)$ be random variables defined on probability space $(\Omega, \CF,\PP)$ with joint distribution $\mu^\e$. Here $\bar q_i^\e :=(\bar \Phi_i^\e, \bar \Psi_i^\e)$ is a point in $\CM \times \CM^2, i=0,\dots, N-1$. Thus $\bar q_i^\e$ has the distribution $\bar \mu_i^\e(dq_i):=[\mu^\e]_{i|0,\dots,i-1}(dq_i|\bar q_0^\e, \dots, \bar q_{i-1}^\e)$ conditioned on the filtration $\bar \CF_i^\e:=\sg\{\bar q_j^\e, j=0,\dots, i-1\}$. Starting at $\bar X_0^\e=x_0$, we define $\bar X_i^\e$ recursively by 
\beq{cont-path}
\bar X_{i+1}^\e =\bar X_i^\e + \e U(\bar X_i^\e, \bar \Phi_i^\e, \bar\Psi_i^\e).
\eeq
We denote by $\bar X^\e(t)$ the piecewise linear interpolation of $\{\bar X_i^\e\}$ with $\bar X^\e(t_i)=\bar X_i^\e$. To keep trace of time dependence for ${
\bar \mu_i^\e}$, we define $
\bar \mu^{\e}(B,D,E) = \int_E \bar \mu^{\e}(B,D|t)dt$,
$\bar \mu^\e(B,D|t)= \bar \mu_i^\e(B,D)$,
and
\beq{la-ga}\barray
\ad\!\!\!\!\!\!\!\!\!\! \la^{\e}(A,C, B, D|t)
=\dl_{(\bar\Phi_i^\e,\bar \Psi_i^\e)}(A, C) \bar \mu_i^{\e}(B, D),\\ 
\ad\!\!\!\!\!\!\!\!\!\! \ga^{\e}( A, C, B, D|t) =\dl_{(\bar \Phi_i^\e,\bar \Psi_i^\e)}(A,C) \rho^1_{\bar X_i^\e}(\bar\Phi_i^\e, B)\rho^2(\bar\Psi_i^\e, D),
\earray\eeq
if $t\in [t_i, t_{i+1})$. We let
\beq{la-ga-int}\barray
\la^\e(A, C, B, D, E)=\int_E \la^\e(A, C, B, D|t)dt, \\
\ga^\e(A, C, B, D, E)=\int_E \ga^\e(A, C, B, D|t)dt.
\earray\eeq 
}
\end{const}

The following proposition gives the tightness of sequences $\{\bar X^\e\},\{\bar \mu^\e\}, \{\la^\e\}$, and $\{\ga^\e\}$ and characterizes their limits. The tightness of $\{\bar \mu^\e\}$, $\{\la^\e\}$, and $\{\ga^\e\}$ follows from the compactness of $\CM$ and $[0,T]$. We note that the underlying limiting system \eqref{barX} is an ordinary differential equation subject to algebraic constraints \eqref{marg}, which is a distinct feature in our results.

\begin{prop}\label{prop:tight}
Consider any sequence of controls $\{\bar \mu_i^\e\}$ as in \constref{const}. 
Then the sequences $\{\bar X^\e\}, \{ \bar \mu^
\e\}$, $\{\la^{\e}\}$, and $\{\ga^{\e}\}$ are tight. Consider a convergent subsequence (still denoted by $\e$) of $(\bar X^\e, \bar \mu^{\e}, \la^{\e}, \ga^{\e})$ converging weakly to $(\bar X, \bar \mu, \la, \ga)$, we have $\bar \mu(dz_1,dz_2, dt)$ can be decomposed as $\bar \mu(dz_1,dz_2|t)dt$, where $\bar \mu(dz_1,dz_2|t)$ is a stochastic kernel on $\CM\times\CM^2$ given $t \in [0,T]$, and w.p.1 for all $t\in [0,T]$, 
\beq{barX}
\bar X(t)= x_0+ \int_0^t \int_{\CM^3} U(\bar X(s), z_1,z_2)\bar \mu(dz_1, dz_2, ds).
\eeq
In addition, 
\beq{lim-la}\barray
\la(A, C, B, D, E) =\int_E \la(A, C, B, D | t)dt,
\earray\eeq
and
\beq{rep-ga}\barray
\ga(A, C, B, D, E)\\
\disp =\int_E \Big(\int_{A\times C} \rho^1_{\bar X(t)}(z_1, B) \rho^2(z_2, D) \bar \mu(dz_1, dz_2|t) \Big) dt.
\earray\eeq
Moreover, w.p.1 for every $t\in [0,T]$, $(z_1,z_2)\in\CM^3$,
\beq{marg} 
\!\!\!\!\bar \mu(z_1,z_2)=\!\!\!\!\!\!\sum_{z_3\in \CM,
z_4\in \CM^2}\!\!\! \rho_{\bar X(t)}^1(z_1,z_3)\rho^2(z_2,z_4)\bar \mu(z_3,z_4).
\eeq
\end{prop}
\para{Proof of \propref{prop:tight}.}
The proof is omitted due to page limitations; we refer the reader to  \cite[Theorem 5.3.5]{DE97}.
\qed

\begin{prop}\label{prop:var-form}
For any bounded continuous function $F: C_{x_0}([0,T]; \rr^{d}) \to \rr$, we have the variational formula:
\beq{var-form}\barray
\ad\!\!\!\!\!\!\! -\e \log \EE e^{-\e^{-1} F(X^\e)}\\
\aad\!\!\!\!\!\!\! =\inf_{\{\bar \mu_i^\e \}} \EE [F(\bar X^\e) +\e \!\! \sum_{i=0}^{T/\e -1} \! R (\bar \mu_i^\e (\cdot,\cdot) ||\rho^1_{\bar X_i^\e}(\bar\Phi_i^\e,\cdot) \rho^2(\bar \Psi_i^\e,\cdot)) ],
\earray\eeq
where $R(\cdot||\cdot)$ is the relative entropy, $\{\bar \mu_i^\e\}$ is a collection of random probability measure satisfying the conditions: $\bar \mu_i^\e$ is measurable with respect to the $\sg$-algebra $
\bar \CF_i^\e=\sg\{\bar\Phi_0^\e, \bar\Psi_0^\e, \cdots, \bar\Phi_{i-1}^\e, \bar \Psi_{i-1}^\e\}
$ and given $\bar \CF_{i}^\e$, the conditional distribution of $(\bar \Phi_i^\e,\bar \Psi_i^\e)$ is $\bar \mu_i^\e$. The $\bar X_i^\e$ is defined via recursion \eqref{sa-z} using the controlled sequences $\{\bar\Phi_i^\e, \bar \Psi_i^\e\}$.
\end{prop}

\para{Proof of \propref{prop:var-form}.}
The proof follows that of \cite[Theorem 4.5]{BD19}, thus the details are omitted.
\qed 

Let us define
\beq{L}\barray
L(x,\beta)\disp
=\inf_{\mu}\bigg\{\inf_{\wdt \ga \in \CA(\mu)} R(\wdt \ga \| \mu \otimes \rho^1_{x}\otimes \rho^2 ): \\ 
\qquad\qquad\qquad\quad \disp 
\beta =\int_{\CM^3} U(x,z_1,z_2)\mu(dz_1,dz_2) \\
\text{ s.t. algebraic constraint: } \forall (z_1,z_2)\in \CM^3, \\
\disp\mu(z_1,z_2)=\sum_{z_3\in \CM,
z_4\in \CM^2} \rho_{x}^1(z_1,z_3)\rho^2(z_2,z_4)\mu(z_3,z_4)
\Big\},
\earray 
\eeq
where $\CA(\mu):=\{\wdt \ga \in \CP(\CM^3 \times \CM^3): [\wdt \ga]_1=[\wdt \ga]_2=\mu\}$. Here $[\wdt \ga]_1$ and $[\wdt \ga]_2$ are  the first and second marginals of $\wdt \ga$. 

{\bf Main Result.}
With the variational formula in \propref{prop:var-form}, we have the following main result.

\begin{thm}\label{thm:ldp-x}
The sequence $\{X^\e(t)\}, t\in [0,T]$ satisfies the large deviations principle with speed $\e^{-1}$ and the rate function $I: C_{x_0}([0,T];\rr^{d}) \to \rr$ having the form 
\beq{I-x}
I(\phi)=\left\{\barray 
\disp \int_0^T L(\phi,\dot \phi)dt, \quad \text{ if } \phi \in \mathcal{AC}_{x_0}([0,T];\rr^{d}),\\
+\infty,\quad \text{otherwise},
\earray \right.
\eeq
where $\mathcal{AC}_{x_0}([0,T];\rr^{d})$ is set of absolutely continuous function on $[0,T]$ with initial value $x_0$.
\end{thm}

\begin{rem}
The function $I$ defined in \eqref{I-x} is indeed a rate function and the map $(x,\beta) \mapsto L(x,\beta)$ is lower semicontinuous on $\rr^d\times \rr^d$, and for each $x\in \rr^d$, the map $\beta\mapsto L(x,\beta)$ is convex; see \cite[Theorem 4.13, Lemma 4.14]{BD19}. 
\end{rem}

\subsection*{Example. Escape Probabilities}
An important application of \thmref{thm:ldp-x} is to compute the escape probability and mean exit time. Let $X^*$ be the asymptotically stable point of $\dot{X}(t)=\bar U(X(t))$ with $\bar G$ in the domain of attraction of $X^*$. For $x\in G$, define $\tau_G^\e=\inf\{t: X^\e(t)\notin G\}$, then we have
$$
\lim_{\e\to 0}\e \log \PP_{x}\{\tau_G^\e \leq T\}=-\inf_{\phi\in \CO} I(\phi),$$
where $\CO=\{\phi\cd: \phi(0)=x, \phi(t)\in \partial G \text{ for some } t \leq T\}$. Moreover, define $I_0(\phi)=\inf\{I(\phi): \phi(0)=x\}$, we get
\bea 
\lim_{\e\to 0}\log\EE_x \tau_G^\e=I_0(x)
\eea

In the setting of stochastic gradient descent (SGD) algorithm in neural network, our result can also be used to show the exponential escape efficiency of SGD from sharp minima for more general noise than additive Gaussian noise considered in \cite{II21}; see also \cite{XSS20, ZWYWM18}. Again, the state space of our noise can be non-compact. For tracking systems, large deviations result can also apply to measure the tracking ability; see \cite[Chapter 10, 11]{Ku84}.

\section{Large deviations upper bound}\label{sec:up-bound}
In this section, we prove the large deviations result in  \thmref{thm:ldp-x} using Laplace principle. We first prove the Laplace principle upper bound, the lower bound requires constructing a sequence of nearly optimal controls, which is involved and needs more care. Thus we leave it to \secref{sec:low-bound}. The idea to prove the Laplace principle upper bound is to use the definition of infimum in the variational representation formula \eqref{var-form} to obtain a sequence of $\wdt \e$-optimal controls and then take liminf using the characterization in \propref{prop:tight}. 

\begin{prop}\label{prop:up}
Under assumptions (A1) and (A2), for any bounded continuous function $F: C_{x_0}([0,T]; \rr^{d})\to \rr$, we have the Laplace principle upper bound
\bea\ad 
\liminf_{\e\to 0} -\e\log \EE e^{- \e^{-1} F(X^\e)} \geq \inf_{\phi} \big\{F(\phi)+I(\phi)\big\}.
\eea
where the infimum is taken over the set $ \mathcal{AC}_{x_0}([0,T];\rr^d)$.
\end{prop}

\para{Proof of \propref{prop:up}.}
From the variational representation formula in \eqref{var-form}, for any $\wdt \e>0$, the infimum implies there exists a sequence of controls $\{\bar \mu_i^\e \}$ such that
\beq{up-e}\barray\!\!\!\!\!
-\e\log\EE e^{-\e^{-1} F(X^\e)} +\wdt \e \\
\!\!\!\!\disp \geq \EE[ F(\bar X^\e) +\e \sum_{i=0}^{T/\e -1} R(\bar \mu_i^\e (\cdot,\cdot) || \rho^1_{\bar X_i^\e}(\bar \Phi_i^\e, \cdot)\rho^2(\bar\Psi_i^\e,\cdot)) ].
\earray\eeq
From the chain rule for relative entropy in \cite{BD19}, we have 
\beq{CR-1}\barray
R(\bar \mu_i^{\e}(\cdot,\cdot) || \rho^1_{\bar X_i^\e}(\bar \Phi_i^\e,\cdot)\rho^2(\bar\Psi_i^\e,\cdot)) \\
= R(\dl_{(\bar\Phi_i^\e,\bar \Psi_i^\e)}(\cdot,\cdot) || \dl_{(\bar\Phi_i^\e,\bar \Psi_i^\e)}(\cdot,\cdot)) \\
\quad + R(\bar \mu_i^{\e}(\cdot,\cdot)  || \rho_{\bar X_i^\e}^1(\bar\Phi_i^\e,\cdot) \rho^2(\bar\Psi_i^\e,\cdot)) \\
= R(\dl_{(\bar\Phi_i^\e,\bar \Psi_i^\e)}(\cdot,\cdot) \bar \mu_i^\e(\cdot, \cdot) ||  \dl_{(\bar\Phi_i^\e,\bar\Psi_i^\e)}(\cdot, \cdot)  \\
\qquad\qquad\qquad\qquad\qquad\quad \quad \otimes \rho^1_{\bar X_i^\e}(\bar\Phi_i^\e, \cdot) \rho^2(\bar\Psi_i^\e, \cdot)) \\
= R(\la^{\e}(\cdot, \cdot, \cdot, \cdot|t) || \ga^{\e}(\cdot, \cdot, \cdot, \cdot|t)).
\earray\eeq
By the tightness of $\{\la^\e\},\{\ga^\e\}$ in \propref{prop:tight}, there is a convergent subsequence converging to $\{\la\}, \{\ga\}$, respectively. We work with such convergent subsequence and still denote it by $\e$ for notation simplicity. Thus, taking liminf in \eqref{up-e} gives the following 
\beq{lim-upe}\barray
\ad
\liminf_{\e\to 0} -\e \EE e^{-\e^{-1} F(X^\e)} +\wdt \e \\
\aad \geq \liminf_{\e \to 0}\EE [F(\bar X^\e)+\sum_{i=0}^{T/\e -1} \e R(\la^{\e}(\cdot|t)||\ga^{\e}(\cdot|t))] \\
& \geq \EE[F(\bar X)+ \int_0^T R(\la(\cdot|t)||\ga(\cdot|t)dt],
\earray\eeq
where the last inequality is due to the Fatou's lemma and the lower semi-continuity of relative entropy. By \eqref{rep-ga}, 
\beq{rewrite-CR}\barray\disp
\int_0^T R(\la(\cdot|t) \| \ga(\cdot|t))dt \\
=\disp \int_0^T R(\la(dz_1, dz_2,dz_3,dz_4|t) ||  \ga(dz_1, dz_2,dz_3,dz_4|t))dt\\
=\disp \int_0^T R (\la(dz_1, dz_2, dz_3, dz_4|t) || \rho^1_{\bar X(t)}(z_1, dz_3)\\
\qquad\qquad \qquad \qquad\qquad\qquad \otimes \rho^2(z_2,dz_4) \bar \mu(dz_1,dz_2|t) ) dt \\
\geq\disp \int_0^T  L(\bar X(t), \dot{\bar X}(t)) dt.
\earray\eeq
Combining \eqref{lim-upe} and \eqref{rewrite-CR} yields
\bea \ad 
\liminf_{\e \to 0} \big[ -\e \log \EE e^{-\e^{-1} F(X^\e)} \big]+\wdt \e \\
& \geq \EE [F(\bar X(t))+\int_0^T  L(\bar X(t),\dot{\bar X}(t)) dt ] \\
& \geq \inf_{\phi} \{F(\phi)+\int_0^T  L(\phi,\dot{\phi}(t)) dt \},
\eea
where the infimum is taken over $ \mathcal{AC}_{x_0}([0,T];\rr^{d})$. Since $\wdt \e$ is arbitrarily small, the Laplace upper bound is proved.
\qed

\section{Large deviations lower bound}\label{sec:low-bound}
\begin{prop}\label{prop:lower}
Under our assumptions (A1) and (A2), for any bounded and continuous function $F: C_{x_0}([0,T];\rr^{d}) \to \rr$, we have the Laplace principle lower bound:
\beq{lower-bd1}
\limsup_{\e \to 0}- \e \log\EE[ e^{- \e^{-1} F(X^\e)}] \leq \inf_{\phi} \{F(\phi)+I(\phi)\},
\eeq
where the infimum is taken over $\mathcal{AC}_{x_0}([0,T]; \rr^d)$. 
\end{prop}

The idea to prove the Laplace principle lower bound \eqref{lower-bd1} is to construct a nearly optimal trajectory $\bar \phi \in \mathcal{AC}_{x_0}([0,T];\rr^d)$. Based on $\bar \phi$, we construct a sequence of near optimal controls to be used in the representation form \eqref{var-form} such that the running cost (the term involving relative entropy) is close to $I(\bar \phi)$ and the corresponding controlled process converges to the nearly optimal trajectory $\bar \phi$. To this purpose, we assume $\inf_{\phi} \{F(\phi)+I(\phi)\}<\infty$, otherwise \eqref{lower-bd1} is trivial.

Fix a small parameter $\wdt \e>0$, the infimum in $\inf_\phi \{F(\phi)+I(\phi)\}$ implies there exists a $\zeta \in \mathcal{AC}_{x_0}([0,T]; \rr^d)$ such that
\bea \ad 
F(\zeta)+I(\zeta)\leq \inf_{\phi\in \mathcal{AC}_{x_0}([0,T];\rr^d)}\{F(\phi)+I(\phi)\} + \wdt \e.
\eea
We note that the function $\zeta\cd$ is bounded due to its continuity on finite interval $[0,T]$. In fact, without loss of generality, we can also assume the derivative $\dot{\zeta}\cd$ is bounded. We omit the argument of this assumption but refer to \cite[Lemma 4.17]{BD19}. Thus, there exist $M_1<\infty$ and $M_2<\infty$ such that
\bea\ad
\sup_{t\in [0,T]}\|\zeta(t)\| \vee \sup_{t\in [0,T]} \|\dot{\zeta}(t)\| \leq M_1 <\infty,
\eea
and 
\bea\ad 
\sup_{\{(x,\beta): \|x\| \leq M_1+1, \|\beta\|\leq M_1+1 \} } L(x,\beta) \leq M_2 <\infty.
\eea
For  any $\bar \dl >0$, we define $\zeta^{\bar \dl}$ as the piecewise linear interpolation of $\zeta$. Then the absolutely continuity of $\zeta$ implies that $\dot{\zeta}^{\bar\dl}(t)$ converges to $\dot{\zeta}(t)$ for a.e. $t\in [0,T]$, as $\bar{\dl} \to 0$. By the continuity of $L$ and the dominated convergence theorem, we conclude that there is a $\bar{\dl}_1>0$ such that $
F(\zeta^{\bar{\dl}_1})+I(\zeta^{\bar{\dl}_1})\leq F(\zeta)+I(\zeta)+\wdt \e.$
We then take $\bar\phi=\zeta^{\bar{\dl}_1}$. Thus, 
$
F(\bar\phi)+I(\bar\phi)\leq F(\zeta)+I(\zeta)+\wdt \e.
$
Before we construct the control based on $\bar\phi$, we need the following lemma. 

\begin{lem}\label{lem:inf-L}
Suppose (A1) and (A2) hold. For $\wdt \e>0$, let $(x,\beta)\in \rr^{d}\times \rr^{d}$ such that $L(x,\beta)<\infty$, then there exists a probability measure $\nu^{x,\beta}$ such that
\beq{inf-R}
\inf_{\wdt \ga \in \CA(\nu^{x,\beta})}\ R(\wdt \ga || \nu^{x,\beta} \otimes \rho^1_x \otimes \rho^2) \leq L(x,\beta)+\wdt \e,
\eeq
and
\beq{nu-U}\barray
\beta=\int_{\CM^3} U(x,z_1,z_2) \nu^{x,\beta}(dz_1,dz_2).
\earray\eeq
For any $\dl>0$, define the measure $
\nu^{x,\beta,\dl}:=(1-\dl/2) \nu^{x,\beta} + \dl/2 ( \pi_{x}^1  \otimes \pi^2).
$
Then there exists a transition kernel $p^{x,\beta,\dl}(z_1,z_2,dz_1,dz_2)$ such that $\nu^{x,\beta,\dl}$ is the unique invariant measure of $p^{x,\beta,\dl}$ and the associated Markov chain with respect to $p^{x,\beta,\dl}$ is ergodic. Moreover, 
\beq{R-dl} \barray\ad 
R(\nu^{x,\beta,\dl} \otimes p^{x,\beta,\dl} || \nu^{x,\beta,\dl} \otimes \rho^1_{x}\otimes \rho^2)\\
\aad \leq \inf_{\wdt \ga\in \CA(\nu^{x,\beta})} R(\wdt \ga|| \nu^{x,\beta} \otimes \rho_{x}^1\otimes \rho^2) \leq L(x,\beta)+\wdt \e.
\earray\eeq
\end{lem}

\para{Proof of \lemref{lem:inf-L}.} The proof is very similar to that of \cite[Lemma 6.17]{BD19} and \cite[Lemma 8.6.3]{DE97}, thus details are omitted. The reason of introducing $\nu^{x,\beta,\dl}$ rather than using $\nu^{x,\beta}$ directly is to guarantee the associated chain is ergodic.

\subsection{Construction of  a near optimal control}
From $p^{x,\beta,\dl}$ in  \lemref{lem:inf-L}, we now in the position to construct the near optimal controls. To proceed, we divide the interval $[0,T]$ into sub-intervals with length $\Dl>0$ and let $T/\Dl$ be an integer for convenience. Define $\tau_k: = k\Dl, k=0,\dots, T/\Dl$. For a given $\dl>0$ and $k=0,\dots, T/\Dl-1$, we define
\bea\ad \!\!\!\!\!\!\!
\nu_i^\e(dz_1,dz_2) \\
\aad\!\!\!\!\!\!\! =\left\{\barray
\rho_{\bar \phi(\tau_k)}^1(\bar\Phi_{i}^\e, dz_1) \rho^2(\bar \Psi_i^\e, dz_2),\; k \Dl/\e  \leq i < k \Dl/\e + \ell_0 \\
p^{\bar \phi(\tau_k), \dot{\bar\phi}(\tau_k), \dl}(\bar \Phi_i^\e,\bar \Psi_i^\e, dz_1,dz_2),\\
\qquad \qquad\qquad\quad k\Dl/\e+\ell_0 \leq i \leq (k+1)\Dl/\e -1,
\earray\right.
\eea
where $i=k\Dl/\e, \dots, (k+1)\Dl/\e-1$ and $\ell_0$ is defined in \eqref{tr-cond}. Starting from $(\bar \Phi_0^\e,\bar\Psi_0^\e)$, this $\nu_i^\e$ defines a controlled sequences of $(\bar\Phi_i^\e, \bar\Psi_i^\e)$ (abusing the notation) such that the conditional distribution of $(\bar\Phi_i^\e, \bar\Psi_i^\e)$ conditioned on $\bar \CF_{i}^\e=\sg\{\bar\Phi_j^\e, \bar\Psi_j^\e, j=0,\dots, i-1\}$ is $\nu_i^\e$. The corresponding controlled state process via recursion \eqref{sa-z} using controlled sequence $\{\bar \Phi_i^\e, \bar \Psi_i^\e\}$ is denoted by $\bar X_i^\e$ (abusing notation). To make sure the controlled state process $\bar X_i^\e$ will not be far away from $\bar\phi$, we define 
the stopping time $\bar \tau^\e$ as
\beq{bar-tau}
\bar \tau^\e:=\inf\{i\in \NN: \|\bar X_i^\e -\bar \phi(t_i)\| > 1\} \wedge (T/\e).
\eeq
We set 
\bea
\bar \nu_i^\e(dz_1, dz_2)=\left\{\barray
\nu_i^\e (dz_1,dz_2), \quad   \ad  i <  \bar \tau^\e, \\
\rho_{\bar X_{i}^\e}(\bar\Phi_{i}^\e,dz_1)\rho^2(\bar\Psi_{i}^\e, dz_2),\quad  \ad  \text{otherwise}.
\earray\right.
\eea
For the rest of this section, we will work with the sequence of controls $\{\bar\nu_i^\e\}$. Likewise, to keep record of the time, we define $\bar \nu^\e(B, D|t):= \bar \nu_i^\e(B,D)$ if $t\in [t_i, t_{i+1})$ and $\bar \nu^\e \in \CP(\CM^3\times[0,T])$ by
$
\bar \nu^\e(B, D, E)=\int_E \bar \nu^\e(B, D|t)dt.
$
\begin{rem}
The sequences $\{\bar\nu^\e\}$ is tight due to the compactness of $\CP(\CM^3)$. Moreover, $[0,T]$ is a compact interval, thus $\{\bar\tau^\e\}$ is also tight. We will work with a convergent subsequence of $\{\bar\tau^\e\}$ (still denote by $\e$) and let $\bar\tau$ be its limit.
\end{rem}

\subsection{Convergence of controls and the controlled processes} 
For each $t\in [0,T]$, we denote by $\nu^{\bar \phi(t),\dot{\bar\phi}(t),\dl}$ the associated measure in \lemref{lem:inf-L} such that \eqref{R-dl} holds. The corresponding transition kernel is denoted by $p^{\bar\phi(t),\dot{\bar\phi}(t),\dl}$. Let us define
\beq{def-nu}\barray
 \nu(B, D, E)=\int_E  \nu(B, D|t)dt,
\earray\eeq
where $
\nu(B, D|t)=  \nu^{\bar\phi(\tau_i),\dot{\bar\phi}(\tau_i),\dl}(B, D),\ t\in [\tau_i, \tau_{i+1}),\ t \leq \bar\tau$. The following lemma establishes the weak limit of $\bar\nu^\e$ and that of the controlled process $\bar X^\e$ using $\bar\nu^\e$. Here $\bar X^\e$ is the piecewise linear interpolation of $\{\bar X_i^\e\}$.

\begin{lem}\label{lem:conv-control}
Suppose assumptions (A1) and (A2) hold, then for every sequence of $(\bar \nu^\e, \bar X^\e)$, there exists a further subsequence that converges weakly to $(\nu,\bar \phi)$, where $\nu$ is defined in \eqref{def-nu} and from \eqref{nu-U}, $\bar\phi$ has the following form:
\beq{b-phi}\barray
\bar\phi(t) & =\bar\phi(0) \\
& +\int_0^t \int_{\CM^3}  U(\bar\phi(s),z_1,z_2) \nu^{\bar\phi(s),\dot{\bar\phi}(s)}(dz_1,dz_2|s) ds.
\earray\eeq
\end{lem}

\para{Proof of \lemref{lem:conv-control}.} The proof of the lemma can be found in \cite{BD19, DE97}. But for completeness, we put it in the appendix.

\para{Proof of \propref{prop:lower}.}
By \lemref{lem:conv-control}, we have \allowdisplaybreaks
\bea \ad\!\!\!\!\!\!
\limsup_{\e\to 0} -\e \log\EE e^{-\e^{-1} F(X^\e)} \\
\aad\!\!\!\!\!\!\leq \limsup_{\e \to 0} \EE \big[F(\bar X^\e)\\
\aad\qquad \qquad\quad + \frac{1}{\e}\sum_{i=0}^{T/\e-1} R(\bar \nu_i^\e(\cdot,\cdot)|| \rho_{\bar X_i^\e}(\bar \Phi_i^\e,\cdot) \rho^2(\bar\Psi_i^\e,\cdot)\big] \\
\aad \!\!\!\!\!\! =\limsup_{\e\to 0} \EE \big[F(\bar X^\e)+ \sum_{k=0}^{T/\Dl-1} \Dl \frac{\e}{\Dl}\\
\aad\qquad \times \sum_{i=k\Dl/\e}^{(k+1)\Dl/\e-1} R(\bar\nu_i^\e(\cdot,\cdot) || \rho_{\bar X_i^\e}^1 (\bar\Phi_i^\e, \cdot)\rho^2(\bar\Psi_i^\e,\cdot) )\big] \\
\aad\!\!\!\!\!\! \leq \limsup_{\e\to 0} \EE \big[ F(\bar X^\e)+ \sum_{k=0}^{T/\Dl-1} \Dl \frac{\e}{\Dl} \sum_{i=k\Dl/\e+\ell_0}^{(k+1)\Dl/\e-1} \\
\aad R(p^{\bar\phi(\tau_k),\dot{\bar\phi}(\tau_k),\dl}(\bar\Phi_i^\e,\bar\Psi_i^\e,\cdot, \cdot) || \rho_{\bar X_i^\e}(\bar\Phi_i^\e, \cdot) \rho^2(\bar\Psi_i^\e,\cdot))  \big] \\
\aad\!\!\!\!\!\! \leq \limsup_{\e\to 0}\EE \big[F(\bar X^\e) + \Dl \sum_{k=0}^{T/\Dl-1}\frac{\e}{\Dl} \sum_{i=k\Dl/\e+\ell_0}^{(k+1)\Dl/\e-1} \\
\aad  R(p^{\bar\phi(\tau_k),\dot{\bar\phi}(\tau_k),\dl}(\bar\Phi_i^\e,\bar\Psi_i^\e, \cdot, \cdot) || \rho_{\bar\phi(\tau_k)}^1(\bar\Phi_i^\e,\cdot)\rho^2(\bar\Psi_i^\e, \cdot) )\big] \\
\aad\!\!\!\!\!\! \leq \EE [F(\bar x^\Dl)+ \Dl\\
\aad\; \times  \sum_{k=0}^{T/\Dl-1} \int_{\CM^3} R(p^{\bar\phi(\tau_k),\dot{\bar\phi}(\tau_k),\dl}(z_1,z_2,\cdot,\cdot) || \rho_{\bar\phi(\tau_k)}^1(z_1,\cdot) \\
\aad\qquad\qquad\qquad \qquad  \otimes \rho^2(z_2,\cdot))\nu^{\bar\phi(\tau_k),\dot{\bar\phi}(\tau_k),\dl}(dz_1,dz_2)].
\eea
The last inequality is because $\nu^{\bar \phi(\tau_k), \dot{\bar\phi}(\tau_k),\dl}$ is the invariant measure for transition kernel $p^{\bar\phi(\tau_k), \dot{\bar\phi}(\tau_k),\dl}$. The third inequality is because $\bar X_i^\e, k\Dl/\e + \ell_0 \leq i \leq (k+1)\Dl/\e -1$  converges to $\bar\phi(\tau_k)$. 
The second inequality is from the definition of $\bar \nu_i^\e$. 
By the chain rule for relative entropy \cite[Corollary 2.7]{BD19}, 
\bea\barray
&\EE\big[\int_{\CM^3} R(p^{\bar\phi(\tau_k),\dot{\bar\phi}(\tau_k),\dl}(z_1,z_2,\cdot,\cdot) || \rho_{\bar\phi (\tau_k)}(z_1,\cdot)\\
& \qquad\qquad\otimes \rho^2(z_2,\cdot))
\nu^{\bar\phi(\tau_k),\dot{\bar\phi}(\tau_k),\dl}(dz_1,dz_2) \big]\\
& =\EE \big[R(\nu^{\bar\phi(\tau_k),\dot{\bar\phi}(\tau_k),\dl}\otimes p^{\bar\phi(\tau_k),\dot{\bar\phi}(\tau_k),\dl}\\
& \qquad\qquad \qquad || \nu^{\bar\phi(\tau_k),\dot{\bar\phi}(\tau_k),\dl}\otimes \rho_{\bar\phi(\tau_k)}^1 \otimes \rho^2 \big] \\
& \leq \EE L(\bar\phi(\tau_k),\dot{\bar\phi}(\tau_k))+\wdt \e.
\earray\eea
Therefore, we have 
\bea\ad 
\limsup_{\e\to 0} -\e \log\EE e^{-\e^{-1}F(X^\e)} \\
\aad \leq \EE [F(\bar x^\Dl)+ \Dl \sum_{k=0}^{T/\Dl-1} L(\bar\phi(\tau_k),\dot{\bar\phi}(\tau_k))]+\wdt \e.
\eea
Taking $\Dl \to 0$ and $\dl\to 0$, we have
\bea \ad 
\limsup_{\e\to 0}-\e \log\EE e^{-\e^{-1}F(X^\e)} \\
& \leq E [F(\bar\phi)+ \int_0^T L(\bar\phi(s),\dot{\bar\phi}(s))ds ] +\wdt \e \\
\aad =  F(\bar\phi)+I(\bar\phi)+\wdt \e \leq \inf_{\phi}(F(\phi)+I(\phi))+2\wdt \e.
\eea
Since $\wdt \e$ is arbitrarily small, we obtain
\bea\ad 
\limsup_{\e\to 0}-\e \log\EE e^{-\e^{-1} F(X^\e)} \leq \inf_{\phi}(F(\phi)+I(\phi)).
\eea
The proof of large deviations lower bound is completed. 
\qed

\section{Conclusions and Discussion}
The main result of the paper was \thmref{thm:ldp-x} which gives the large deviations result for the regret minimization stochastic approximation algorithm. Recall that the weak convergent limit of the algorithm is an algebraically constrained ordinary differential equation. We see that the rate function in \eqref{L} depends on this algebraic constraint. As the example below Theorem~\ref{thm:ldp-x}, an important application of the large deviations analysis is to compute the escape probability and mean exit time of the stochastic approximation algorithm. Finally, we emphasize that 
the large deviations approach used here is based on weak convergence.

\bibliographystyle{abbrv}

\section{Appendix}
\para{Proof of \lemref{lem:conv-control}.}
The tightness of $\{\bar X^\e \}$ follows from \propref{prop:tight} using the control $\{\bar \nu^\e\}$. We will first show as $\e \to 0$, the convergent subsequence of $\{\bar X^\e\}$ converges to $\bar x^\Dl$, as $\e \to 0$, where 
\beq{barxl}\barray
\bar x^\Dl(t)
& := x_0 + \sum_{i=0}^{k-1} \int_{\tau_i}^{\tau_{i+1}}\int_{\CM^3}  U(\bar x^\Dl(s), z_1, z_2) \\
& \qquad\qquad\qquad \quad \times \nu^{\bar \phi(\tau_i), \dot{\bar \phi}(\tau_i),\dl}(dz_1,dz_2|s)ds \\
& +\int_{\tau_k}^{t} \int_{\CM^3} U(\bar x^\Dl(s), z_1,z_2) \\
& \qquad\qquad\qquad\times  \nu^{\bar \phi(\tau_k),\dot{\bar \phi}(\tau_k),\dl}(dz_1,dz_2|s)ds,
\earray\eeq
if $t\in [\tau_k, \tau_{k+1})$. Then we prove $\bar x^\Dl$ converges to $\bar\phi$ as $\Dl\to 0$ and $\dl\to 0$. To this purpose, we first show $\bar\nu^\e$ converges weakly to $\nu$. It is sufficient to prove that as $\e\to 0$, $
\int_0^T \int_{\CM^3} f(t,z_1,z_2)\bar \nu^\e (dz_1,dz_2,dt)
\to \int_0^T \int_{\CM^3} f(t,z_1,z_2)  \nu(dz_1,dz_2,dt),
$
for any bounded and uniformly continuous function $f: [0,T]\times \CM\times \CM^2 \to \rr$. We set 
$
\nu(B, D|t)=\lim_{\e\to 0} \bar \nu^\e(B,D|t),\; t\geq \bar\tau.
$
The definition of $\nu$ for $t\geq \bar\tau$ implies that we only need to prove the case of $T=\bar\tau$.
From the definition of $\bar \nu^\e$, we note that
$
\int_0^T \int_{\CM^3} f(t,z_1,z_2) \bar \nu^\e(dz_1,dz_2,dt) 
=\sum_{k=0}^{\frac{T}{\Dl}-1} \sum_{i= k\Dl/\e}^{(k+1) \Dl/\e -1} \int_{t_i}^{t_{i+1}} \int_{\CM^3} f(t,z_1,z_2) \bar \nu^\e(dz_1,dz_2,dt).
$
For each $k=0,\dots,T/\Dl-1$, we will prove 
\beq{conv-k}\barray
\!\!\! \Ga_k & \!\!\!\!\! :=\!\!\! \sum_{i=k \Dl/\e}^{(k+1)\Dl/\e -1} \int_{t_i}^{t_{i+1}} \int_{\CM^3} f(t,z_1,z_2) \bar \nu^\e (dz_1,dz_2,dt) \\
& \quad \to \int_{\tau_k}^{\tau_{k+1}}\int_{\CM^3} f(t,z_1,z_2)  \nu(dz_1,dz_2,dt),
\earray\eeq
as $\e \to 0$. The definition of $\bar \nu^\e$ in \eqref{def-nu} gives us 
\bea
\Ga_k 
& \!\!\!\!=\sum_{i=k \Dl/\e}^{k \Dl/\e +\ell_0-1} \int_{t_i}^{t_{i+1}} \int_{\CM^3} f(t,z_1,z_2)\rho^1_{\bar \phi(\tau_k)}(\bar \Phi_i^\e, dz_1) \\
& \qquad\qquad\qquad\qquad\qquad\qquad \qquad\quad \otimes \rho^2(\bar \Psi_i^\e, dz_2)dt \\
& \quad + \sum_{i=k\Dl/\e+\ell_0}^{(k+1)\Dl/\e-1} \int_{t_i}^{t_{i+1}} \int_{\CM^3} f(t,z_1,z_2)\\
& \qquad\qquad\qquad\quad \times  p^{\bar \phi(\tau_k), \dot{\bar \phi}(\tau_k), \dl}(\bar \Phi_i^\e, \bar \Psi_i^\e, dz_1,dz_2) dt.
\eea
Let us define
\bea
\!\! \Ga_{k,1} &\!\!\! := \sum_{i=k\Dl/\e}^{k \Dl/\e +\ell_0-1} \int_{t_i}^{t_{i+1}}\int_{\CM^3} f(t,z_1,z_2)\rho^1_{\bar \phi(\tau_k)}(\bar \Phi_i^\e, dz_1) \\
& \qquad\qquad\qquad\qquad\qquad \qquad\qquad \otimes \rho^2(\bar \Psi_i^\e, dz_2)dt, \\
\!\! \Ga_{k,2} & \!\!\!:= \sum_{i=k\Dl/\e+\ell_0}^{(k+1)\Dl/\e-1} \int_{t_i}^{t_{i+1}} \int_{\CM^3} f(t,z_1,z_2) \\
\aad \qquad\qquad\qquad \times p^{\bar \phi(\tau_k), \dot{\bar \phi}(\tau_k), \dl}(\bar \Phi_i^\e,\bar \Psi_i^\e, dz_1, dz_2) dt \\
& \quad  - \sum_{i=k\Dl/\e+\ell_0}^{(k+1)\Dl/\e-1} \int_{t_i}^{t_{i+1}} \int_{\CM^3} f(t_i, z_1,z_2)  \\
& \qquad\qquad\qquad \times p^{\bar \phi(\tau_k), \dot{\bar \phi}(\tau_k), \dl}(\bar \Phi_i^\e,\bar \Psi_i^\e, dz_1,dz_2) dt,\\
\!\!\Ga_{k,3} & \!\!\!:= \sum_{i=k\Dl/\e+\ell_0}^{(k+1)\Dl/\e-1} \int_{t_i}^{t_{i+1}} \int_{\CM^3} f(t_i, z_1,z_2) \\
& \qquad\qquad\qquad \times p^{\bar \phi(\tau_k), \dot{\bar \phi}(\tau_k),\dl}(\bar \Phi_i^\e, \bar \Psi_i^\e, dz_1,dz_2) dt \\ 
& - \sum_{i=k\Dl/\e+\ell_0}^{(k+1)\Dl/\e-1} \e  \int_{\CM^3} \int_{\CM^3} f(t_i, z_1,z_2)  \\
\aad \times p^{\bar \phi(\tau_k), \dot{\bar \phi}(\tau_k), \dl}(\theta, \psi, dz_1,dz_2)  \nu^{\bar \phi(\tau_k),\dot{\bar \phi}(\tau_k),\dl}(d\theta, d\psi),
\eea
\bea 
\!\! \Ga_{k,4} \ad\!\!\! :=\sum_{i=k \Dl/\e+\ell_0}^{(k+1) \Dl/\e-1} \e  \int_{\CM^3}f(t_i, \theta,\psi)  \nu^{\bar \phi(\tau_k),\dot{\bar \phi}(\tau_k),\dl}(d\theta, d\psi) \\
\aad - \int_{\tau_k}^{\tau_{k+1}}\int_{\CM^3} f(t,z_1,z_2)  \nu^{\bar\phi(\tau_k), \dot{\bar\phi}(\tau_k),\dl}(dz_1,dz_2,dt).
\eea
It follows that
\beq{Gak-dec}\barray
\Ga_k & = \sum_{i=1}^4 \Ga_{k,i}\\
& +\int_{\tau_k}^{\tau_{k+1}} \int_{\CM^3} f(t,z_1,z_2) \nu^{\bar\phi(\tau_k),\dot{\bar\phi}(\tau_k),\dl} (dz_1, dz_2, dt). 
\earray\eeq
In the above decomposition, we used the result in \lemref{lem:inf-L} that the measure $ \nu^{\bar\phi(t),\dot{\bar \phi}(t),\dl}$ is the invariant measure for the transition kernel $p^{\bar\phi(t),\dot{\bar\phi}(t),\dl}$ for $t\in [0,T]$. It gives 
\bea\barray
& \int_{\CM^3}\int_{\CM^3} f(t_i, z_1,z_2) p^{\bar\phi(\tau_k), \dot{\bar\phi}(\tau_k),\dl}(\theta,\psi, dz_1,dz_2)\\
& \qquad\qquad\qquad\qquad \qquad\quad  \otimes  \nu^{\bar\phi(\tau_k),\dot{\bar\phi}(\tau_k),\dl}(d\theta, d\psi) \\
& = \int_{\CM^3} f(t_i,\theta,\psi) \nu^{\bar\phi(\tau_k),\dot{\bar\phi}(\tau_k),\dl}(d\theta, d\psi).
\earray\eea
Thus, \eqref{Gak-dec} is valid. To show the convergence result \eqref{conv-k}, it is sufficient to prove $\Ga_{k,i}\to 0, i=1,2,3,4$ as $\e \to 0$. For $\Ga_{k,1}$ and $\Ga_{k,2}$, the boundedness of function $f$ and its uniform continuity implies $\Ga_{k,1}\to 0$ and $\Ga_{k,2}\to 0$, as $\e\to 0$. For $\Ga_{k,3}$, we note that 
\bea\barray & 
\e \int_{\CM^3} f(t_i, z_1, z_2) p^{\bar\phi(\tau_k), \dot{\bar\phi}(\tau_k), \dl}(\bar\Phi_i^\e, \bar \Psi_i^\e, dz_1, dz_2)\\
& -\e\int_{\CM^3}\int_{\CM^3} f(t_i, z_1, z_2)p^{\bar\phi(\tau_k), \dot{\bar\phi}(\tau_k),\dl}(\phi, \psi, dz_1, dz_2)\\
& \qquad\qquad\qquad \qquad\qquad\qquad \otimes  \nu^{\bar\phi(\tau_k), \dot{\bar\phi},\dl}(d\theta, d\psi)
\earray\eea
forms a martingale difference sequences. Then, the strong law of large numbers implies $\Ga_{k,3} \to 0$ w.p.1 as $\e\to 0$. For $\Ga_{k,4}$, we note that the function
$
t \mapsto \int_{\CM^3} f(t,z_1,z_2) \nu^{\bar\phi(\tau_k),\dot{\bar\phi}(\tau_k),\dl}(dz_1,dz_2|t)
$
is Riemann integrable. Thus, we view
\bea\barray
\sum_{i=k\Dl/\e+\ell_0}^{(k+1)\Dl/\e -1}\e \int_{\CM^3} f(t_i,\theta,\psi)  \nu^{\bar\phi(\tau_k), \dot{\bar\phi}(\tau_k),\dl)}(d\theta,d\psi)
\earray\eea
as the Riemann sum on the interval $[k\Dl, (k+1)\Dl]$ with step size $\e$. It gives that $\|\Ga_{k,4}\| \to 0$, as $\e\to 0$. These end up with the convergence of $\bar\nu^\e$ to $\nu$. Then, similar to \propref{prop:tight}, it implies that $\bar X^\e$ converges to $\bar x^\Dl$ as $\e\to 0$. The next step is to show $\bar x^\Dl$ converges to $\bar\phi$ as $\Dl\to 0$ and $\dl\to 0$. Define 
\bea
\La(t) &\!\!\!\!\!:= \sum_{i=0}^{k-1} \int_{\tau_i}^{\tau_{i+1}}\int_{\CM^3} U(\bar \phi(s),z_1,z_2)\\
&\qquad\qquad\qquad\qquad \times  \nu^{\bar\phi(\tau_i), \dot{\bar\phi}(\tau_i),\dl}(dz_1,dz_2,ds) \\
&\!\!\! +\int_{\tau_k}^t \int_{\CM^3} U(\bar\phi(s),z_1, z_2)  \nu^{\bar\phi(\tau_k),\dot{\bar\phi}(\tau_k),\dl}(dz_1,dz_2,ds) \\
&\!\!\! - \int_0^t \int_{\CM^3} U(\bar\phi(s),z_1,z_2) \nu^{\bar\phi(s),\dot{\bar\phi}(s)}(dz_1,dz_2, ds).
\eea
From \eqref{barxl} and \eqref{b-phi}, the uniform Lipschitz continuity of $U$ yields
\bea &
\|\bar x^\Dl(t)-\bar \phi(t)\| \\
& \leq \sum_{i=0}^{k-1}\int_{\tau_i}^{\tau_{i+1}}\int_{\CM^3} \|U(\bar x^\Dl(s),z_1,z_2) -U(\bar\phi(s),z_1,z_2)\| \\
&\qquad\qquad\qquad\qquad\qquad  \times  \nu^{\bar\phi(\tau_i),\dot{\bar\phi}(\tau_i),\dl}(dz_1,dz_2|s)ds  \\ 
&  +\int_{\tau_k}^t \int_{\CM^3} \|U(\bar x^\Dl(s),z_1,z_2) -U(\bar\phi(s),z_1, z_2)\|  \\
\aad\qquad\qquad\times \nu^{\bar\phi(\tau_k), \dot{\bar\phi}(\tau_k),\dl}(dz_1,dz_2|s)ds + \|\La(t)\| \\
\aad \leq K \int_0^t \|\bar x^\Dl(s)-\bar\phi(s)\| ds + \|\La(t)\|.
\eea
Thus, the Gronwall inequality implies\bea 
\sup_{t\in [0,T]}\|\bar x^\Dl(t)-\bar \phi(t)\| \leq e^{KT} \sup_{t\in [0,T]}\|\La(t)\|.\eea
Therefore, to prove the convergence of $\bar x^\Dl$ to $\bar\phi$, it is sufficient to prove $\sup_{t\in [0,T]}\|\La(t)\| \to 0$, as $\Dl\to 0$ and $\dl \to 0$. To proceed, we note that 
\bea
& \!\!\!\!\!\!\! \|\La(t)\|\\ 
& \!\!\!\!\!\! \leq \sum_{i=0}^{k-1} \big\| \int_{\tau_i}^{\tau_{i+1}}\int_{\CM^3} U(\bar\phi(s),z_1,z_2) \\
& \!\!\!\!\!\! \times [ \nu^{\bar\phi (\tau_i), \dot{\bar\phi}(\tau_i),\dl}(dz_1,dz_2|s)-  \nu^{\bar\phi (\tau_i),\dot{\bar\phi}(\tau_i)}(dz_1,dz_2|s)]ds \big\| \\
\aad\!\!\!\!\!\! + \big\|\int_{\tau_k}^t \int_{\CM^3} U(\bar \phi(s), z_1, z_2)\\
\aad\!\!\!\! \times \big[ \nu^{\bar\phi(\tau_k),\dot{\bar\phi}(\tau_k),\dl}(dz_1,dz_2|s) - \nu^{\bar\phi (\tau_k),\dot{\bar\phi}(\tau_k)}(dz_1,dz_2|s)\big] ds  \big\| \\
\ad\!\!\!\!\! + \big\|\sum_{i=0}^{k-1}\int_{\tau_i}^{\tau_{i+1}} \int_{\CM^3} U(\bar\phi(s),z_1,z_2)  \nu^{\bar\phi(\tau_i),\dot{\bar\phi}(\tau_i)}(dz_1,dz_2,ds) \\
\aad + \int_{\tau_k}^t \int_{\CM^3} U(\bar\phi(s),z_1,z_2)
\ \nu^{\bar\phi(\tau_k),\dot{\bar\phi}(\tau_k)}(dz_1, dz_2,  ds) \\
\aad- \int_0^t \int_{\CM^3} U(\bar\phi(s),z_1,z_2) \nu^{\bar\phi(s),\dot{\bar\phi}(s)}(dz_1,dz_2,ds) \big\| \\
\ad\!\!\!\!\! =: \sum_{i=0}^{k-1}\|\La_i(t)\|+\|\La_k(t)\|+\|\Theta(t)\|.
\eea
For $\La_i(t)$, we have 
\beq{La-i}\barray\ad\!\!\!\!\!\!\! 
\big\|\int_{\tau_i}^{\tau_{i+1}}\int_{\CM^3} U(\bar \phi(s),z_1,z_2) \big[ \nu^{\bar\phi(\tau_i),\dot{\bar\phi}(\tau_i),\dl}(dz_1,dz_2|s)\\
\aad \qquad\qquad\qquad \qquad\quad\quad -  \nu^{\bar\phi(\tau_i), \dot{\bar\phi}(\tau_i)}(dz_1,dz_2|s)\big] ds \big\| \\
\aad\!\!\!\!\!\!\! \leq \big\|\int_{\tau_i}^{\tau_{i+1}}\int_{\CM^3} U(\bar\phi(\tau_i),z_1,z_2)[ \nu^{\bar\phi(\tau_i),\dot{\bar\phi}(\tau_i),\dl}(dz_1,dz_2|s)\\
\aad \qquad \qquad\qquad\qquad\qquad -\nu^{\bar\phi(\tau_i), \dot{\bar\phi}(\tau_i)}(dz_1,dz_2|s)] ds \big\| \\
\aad\!\!\!\!\!\!\! + \big\| \int_{\tau_i}^{\tau_{i+1}}\int_{\CM^3} \big[U(\bar\phi(s),z_1,z_2) - U(\bar \phi(\tau_i), z_1, z_2)\big]\\
\aad\qquad\qquad\qquad\qquad\qquad  \times  \nu^{\bar\phi(\tau_i),\dot{\bar\phi}(\tau_i),\dl}(dz_1, dz_2|s)ds\big\| \\
\aad \!\!\!\!\!\!\!  + \big\|\int_{\tau_i}^{\tau_{i+1}} \int_{\CM^3} \big[U(\bar\phi(s),z_1,z_2) - U(\bar \phi(\tau_i), z_1, z_2)\big] \\
\aad\qquad\qquad\qquad\qquad\quad \qquad  \times \nu^{\bar\phi(\tau_i),\dot{\bar\phi}(\tau_i)}(dz_1, dz_2|s)ds \big\| \\
\aad\!\!\!\!\!\!\! =: \|\La_{i,1}(t)\| + \|\La_{i,2}(t)\| + \|\La_{i,3}(t)\|.
\earray \eeq
Recall the definition of $\nu^{x,\beta,\dl}$ in \lemref{lem:inf-L}, for $\La_{i,1}(t)$, we have  
\bea\barray &
\|\La_{i,1}(t)\|\\
& =(\tau_{i+1}-\tau_i)\frac{\dl}{2} \\
& \quad  \times \big\|\int_{\CM^3} U(\bar\phi(\tau_i), z_1,z_2) \pi_{\bar\phi(\tau_i)}^1(dz_1)\pi^2(dz_2)  \\
& \qquad\quad  -\int_{\CM^3} U(\bar\phi(\tau_i),z_1,z_2) \nu^{\bar\phi(\tau_i),\dot{\bar\phi}(\tau_i)}(dz_1,dz_2)\big\| \\
& = \Dl \frac{\dl}{2} \big\|\int_{\CM^3} U(\bar\phi(\tau_i), z_1,z_2) \\
& \qquad\qquad\qquad \qquad \times \pi_{\bar\phi(\tau_i)}^1(dz_1)\pi^2(dz_2) - \dot{\bar\phi}(\tau_i) \big\|,
\earray\eea
where the last equality is due to \eqref{nu-U}. For the term $\La_{i,2}$ and $\La_{i,3}$, the uniform Lipschitz continuity of $U$ implies
$
\|\La_{i,2}(t)\| \leq K \int_{\tau_i}^{\tau_{i+1}} \|\bar\phi(s)-\bar\phi(\tau_i)\| ds,
$
and 
$
\|\La_{i,3}(t)\| \leq K \int_{\tau_i}^{\tau_{i+1}} \|\bar\phi(s)-\bar\phi(\tau_i)\|ds.
$
Thus, we obtain for $i=0,\dots, k-1$,
\bea\barray
\|\La_i(t)\|
& \leq \Dl \frac{\dl}{2} \big\|\int_{\CM^3} U(\bar\phi(\tau_i),z_1,z_2)\\
& \qquad\qquad\qquad  \times  \pi^1_{\bar\phi(\tau_i)}(dz_1)\pi^2(dz_2)-\dot{\bar\phi}(\tau_i)\big\| \\
& + 2K \int_{\tau_i}^{\tau_{i+1}} \|\bar\phi(s)-\bar\phi(\tau_i)\| ds. 
\earray\eea
Likewise, we can have the following estimates for $\La_k$,
\bea\barray
\|\La_k(t)\|
& \leq (t-\tau_k)\frac{\dl}{2} \big\|\int_{\CM^3} U(\bar\phi(\tau_k),z_1,z_2)\\
& \qquad\qquad\qquad \times  \pi^1_{\bar\phi(\tau_k)}(dz_1)\pi^2(dz_2) -\dot{\bar\phi}(\tau_k) \big\| \\
& + K \int_{\tau_k}^t \|\bar\phi(s)-\bar\phi(\tau_k)\|ds.
\earray\eea
Now, it is the turn to obtain estimates for $\Theta$. Similar to the argument of \eqref{La-i}, the uniform Lipschitz continuity of $U$ implies that 
\bea
& \!\!\!\!\!\!\!  \|\Theta(t)\| \leq \sum_{i=0}^{k-1} K \int_{\tau_i}^{\tau_{i+1}} \|\bar\phi(s)-\bar\phi(\tau_i)\| ds \\
& \!\!\!\!\!\!\!  + K\int_{\tau_k}^t  \|\bar\phi(s)-\bar\phi(\tau_k)\| ds \\
& \!\!\!\!\!\!\!   + \big\|\sum_{i=0}^{k-1}\int_{\tau_i}^{\tau_{i+1}} \int_{\CM^3} U(\bar\phi(\tau_i),z_1,z_2)\\
& \qquad\qquad\qquad \qquad\qquad\qquad  \times  \nu^{\bar\phi(\tau_i),\dot{\bar\phi}(\tau_i)}(dz_1,dz_2,ds) \\
& \!\!\!\!\!\!\! + \int_{\tau_k}^t \int_{\CM^3} U(\bar\phi(\tau_k),z_1,z_2) \nu^{\bar\phi(\tau_k), \dot{\bar\phi}(\tau_k)}(dz_1, dz_2,ds) \\
& \!\!\!\!\!\!\!  - \int_0^t \int_{\CM^3}  U(\bar\phi(s),z_1,z_2) \nu^{\bar\phi(s),\dot{\bar\phi}(s)}(dz_1, dz_2, ds)\big\| \\
& \!\!\!\!\!\!\!  =: \sum_{i=0}^{k-1} K \int_{\tau_i}^{\tau_{i+1}} \|\bar\phi(s)-\bar\phi(\tau_i)\| ds \\
& + K \int_{\tau_k}^t \|\bar\phi(s)-\bar\phi(\tau_k)\|ds  + \|\Theta_1(t)\|.
\eea
The definition of $\nu^{\bar\phi(\tau_i),\dot{\bar\phi}(\tau_i)}$ gives 
\bea
\disp \|\Theta_1(t)\| = \|\sum_{i=0}^{k-1}(\tau_{i+1}-\tau_i)\dot{\bar\phi}(\tau_i)+(t-\tau_k) \dot{\bar\phi}(\tau_k)-\dot{\bar\phi}(t)\|  
\eea
We note that
$
\sum_{i=0}^{k-1} (\tau_{i+1}-\tau_i) \dot{\bar\phi}(\tau_i) + (t-\tau_k) \dot{\bar\phi}(\tau_k) $
is the piecewise linear interpolation of $\dot{\bar\phi}(t)$. Thus, $\|\Theta_1(t)\| \to 0$ as $\Dl \to 0$. Combining all above estimates, we have
\bea
& \!\!\!\!\!\!\! 
\|\bar x^\Dl - \bar \phi(t)\|_\infty = \sup_{t\in [0,T]} \|\bar x^\Dl(t)-\bar\phi(t)\| \\
& \!\!\!\!\!\!\!  \leq \sum_{i=1}^{T/\Dl-1} \Dl \frac{\dl}{2}\big\| \int_{\CM^3} U(\bar \phi(\tau_i), z_1, z_2)\\
& \qquad\qquad\qquad\qquad \times \pi_{\bar\phi(\tau_i)}^1(dz_1) \pi^2(dz_2)  -\dot{\bar\phi}(\tau_i)\big\| \\
& \!\!\!\!\! +2K \sum_{i=0}^{T/\Dl-1} \int_{\tau_i}^{\tau_{i+1}} \|\bar\phi(s)-\bar\phi(\tau_i)\|ds \\
& \!\!\!\!\!+ \sup_{t\in [0,T]} \|\Theta_1(t)\|.
\eea
For any $\dl>0$, the uniform continuity of $\bar\phi$ implies that for sufficient small $\Dl$, we have 
\bea\ad
\max_{i=0,\dots, T/\Dl-1}\sup_{s\in [\tau_i, \tau_{i+1}]} \|\bar\phi(s)-\bar\phi(\tau_i)\| \leq \dl. 
\eea
It follows that
$
2K \sum_{i=1}^{T/\Dl-1} \int_{\tau_i}^{\tau_{i+1}} \|\bar\phi(s)-\bar \phi(\tau_i)\| ds \leq 2K\dl T \to 0
$
as $\dl\to 0$. Moreover, due to the boundedness of $U$ and $\bar\phi$,  we note as $\dl\to 0$,
\bea\barray &  
\sum_{i=1}^{T/\Dl-1} \Dl \frac{\dl}{2} \big\|\int_{\CM^3} U(\bar\phi(\tau_i), z_1, z_2)\\
& \qquad\qquad\qquad\times  \pi_{\bar\phi(\tau_i)}^1(dz_1) \pi^2(dz_2) -\dot{\bar\phi}(\tau_i)\big\| \to 0.
\earray \eea
Thus, we conclude that $\bar x^\Dl\to \bar\phi $ as $\Dl\to 0$ and $\dl\to 0$.\qed

\end{document}